\documentclass[12pt,a4paper]{article}
\usepackage[english]{babel}
\usepackage{amsmath}
\usepackage{amstext}
\usepackage{amssymb}
\usepackage{amsfonts}
\usepackage{fullpage}
\usepackage{graphicx}
\usepackage{epsf}
\usepackage[dvips]{epsfig}
\usepackage{stmaryrd}
\usepackage[T1]{fontenc}
\usepackage{latexsym}
\usepackage{psfrag}
\usepackage{theorem}

\title{A counter-example to the theorem of Hiemer and Snurnikov}
\author{Thierry Monteil
\footnote{Institut de Math\'ematiques de Luminy, CNRS UPR 9016,
Case 907, 163 Avenue de Luminy, 13288 Marseille cedex 09, France
-- E-Mail: monteil@iml.univ-mrs.fr -- Tel: +33 4 91 26 96 77 --
Fax : +33 4 91 26 96 55 }}
\date{}

\newcommand{\N}{\mathbb N}
\newcommand{\R}{\mathbb R}

\newcommand{\Z}{\mathbb Z}
\newcommand{\Q}{\mathbb Q}

\newenvironment{appli}{\left( \begin{array}{ccc}}{\end{array} \right)}

\newcommand{\ba}{\begin{appli}}
\newcommand{\ea}{\end{appli}}

\newtheorem{thm}{Theorem}
\newtheorem{cor}{Corollary}

\newcommand{\cge}{\xrightarrow[ n \rightarrow \infty]{}}

\renewcommand{\P}{{\mathcal{P}}}

\begin{document}

\maketitle

\begin{abstract}
\noindent A planar polygonal billiard $\P$ is said to have the
finite blocking property if for every pair  $(O,A)$ of points in
$\P$ there exists a finite number of ``blocking'' points $B_1,
\dots , B_n$ such that every billiard trajectory from $O$ to $A$
meets one of the $B_i$'s. As a counter-example to a theorem of
Hiemer and Snurnikov, we construct a family of rational billiards
that lack the finite blocking property.\\

\end{abstract}

\section{Introduction}

\noindent A planar polygonal billiard $\P$ is said to have the
finite blocking property if for every pair  $(O,A)$ of points in
$\P$ there exists a finite number of ``blocking'' points $B_1,
\dots , B_n$ (different from $O$ and $A$) such that every billiard trajectory from $O$ to $A$
meets one of the $B_i$'s.\\

In \cite{HS}, Hiemer and Snurnikov tried to prove that any
rational polygonal billiard has the finite blocking property. The
aim of this paper is to construct a family of rational billiards
that lack the finite blocking property.

\clearpage

\section{The counter-example}\label{example}

\begin{figure}[h]
\begin{center}
\psfrag{1}{$1$}
\psfrag{a}{$\alpha$}
\psfrag{L1}{$L_1$}
\psfrag{L2}{$L_2$}
\psfrag{A(0,2)}{$A(0,2)$}
\psfrag{O}{$O$}
\psfrag{2a}{$2\alpha$}
\psfrag{(a,1-L1)}{$(\alpha,1-L_1)$}
\psfrag{(1,1+L2)}{$(1,1+L_2)$}
\psfrag{(a,1)}{$(\alpha,1)$}
\includegraphics[scale=0.7]{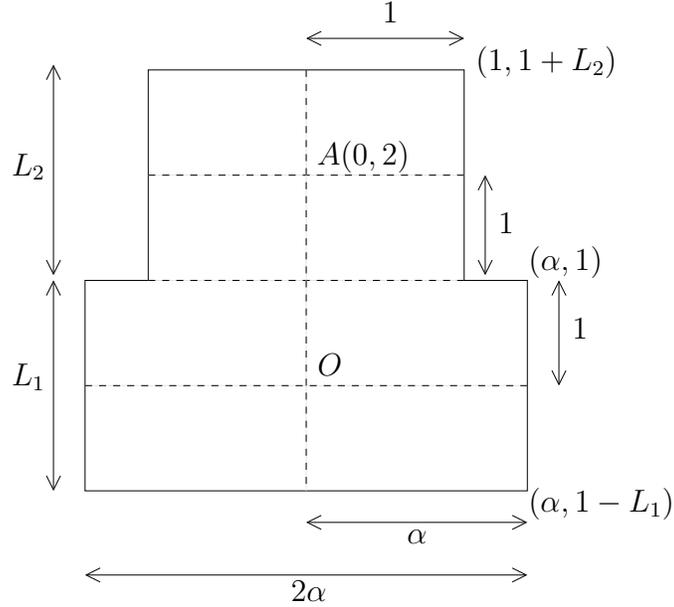}
\caption{\label{def-bill} The polygon $\P_\alpha$.}
\end{center}
\end{figure}

Let $\alpha$ be a positive irrational number 
and $\P_\alpha$ be the polygon drawn in Figure \ref{def-bill}
($L_1$ and $L_2$ can be chosen arbitrarily, greater than 1).

Let $(p_n,q_n)_{n \in \N}$ be a sequence in ${\N^*}^2$ such that:

\begin{itemize}
  \item $q_n$ is strictly increasing
  \item $|p_n - q_n \alpha| < 1$
\end{itemize}

For example, we can take $q_n = n+1$ and $p_n = [ q_n \alpha ]$.

For $n \in \N$, let $\gamma_n$ be the billiard trajectory starting
 from $O$ to $A$ with slope
$$\frac{1}{p_n + q_n \alpha }
=\frac{1}{2q_n \alpha +\lambda_n}=\frac{1}{2p_n -\lambda_n}$$
where $\lambda_n = p_n - q_n \alpha \in ]-1,1[$.

\begin{figure}[!h]

\begin{center}
\psfrag{totalslope}{\small $2(p_n+q_n \alpha)$}
\psfrag{slope1}{\small $2q_n \alpha + \lambda_n$}
\psfrag{slope2}{\small $2p_n - \lambda_n$}
\psfrag{ln}{\small $\lambda_n$}
\includegraphics[width=.98 \linewidth]{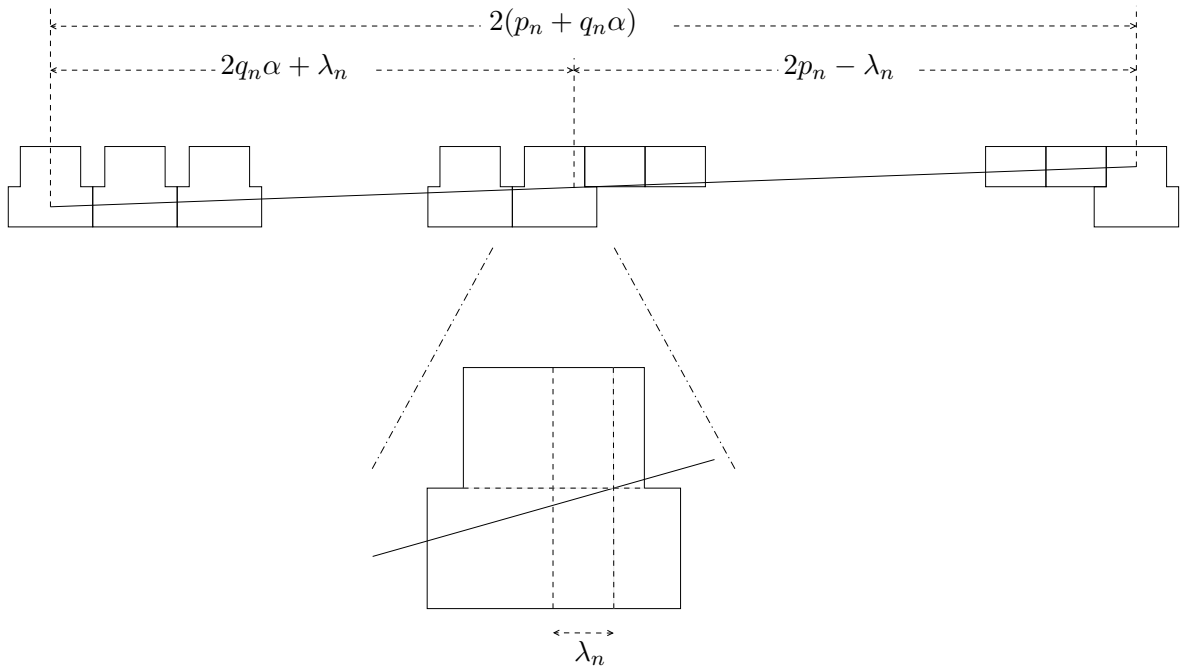}
\caption{\label{unfold} The  unfolding procedure.}
\end{center}
\end{figure}

So, we can check (with the classical unfolding procedure shown in
Figure \ref{unfold}) that
$\gamma_n$ hits $q_n$ walls, passes through $(\lambda_n,1)$, hits $p_n$ walls and then passes through $A(0,2)$.\\

\clearpage

\begin{figure}[!h]
\begin{center}
\includegraphics[scale=0.6]{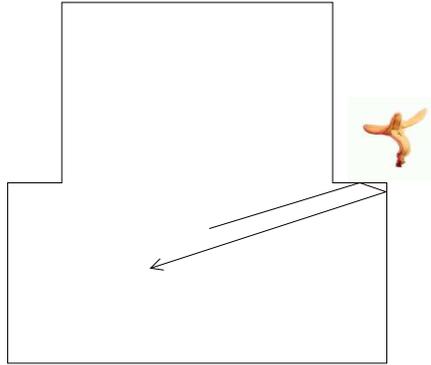}
\caption{\label{banana-peel} The banana peel.}
\end{center}
\end{figure}

The fact that $\lambda_n \in ]-1,1[$ enables us to avoid the banana peel shown in Figure \ref{banana-peel}.\\

Now, we assume by contradiction that there is a point $B(x,y)$ in
$\P_\alpha$ distinct from $O$ and $A$ such that infinitely many
$\gamma_n$ pass through $B$. Hence, there is a subsequence such
that for all $n$ in $\N$, $\gamma_{i_n}$ passes through $B$.

There are two cases to consider:

First case: $y \in ]0,1]$.
By looking at the unfolded version of the trajectory (Figure \ref{unfold}),
we see that $x=\varepsilon_{i_n} y (p_{i_n} + q_{i_n} \alpha ) \ \ [$mod $ 2\alpha]$
where $\varepsilon_{i_n}\in \{-1,1\}$ depends on the parity of the number of bounces of $\gamma_{i_n}$
from $O$ to $B$.

So, there exists a sequence $(k_n)_{n\in \N}$ in $\Z$ such that
$x=\varepsilon_{i_n} y (p_{i_n} + q_{i_n} \alpha) + 2 k_{i_n}
\alpha$.

Taking a further subsequence, we can consider $\varepsilon \circ
i$ to be constant with value $\varepsilon$.

We have $x=\varepsilon y (p_{i_0} + q_{i_0} \alpha) + 2 k_{i_0} \alpha = \varepsilon y (p_{i_1} + q_{i_1} \alpha) + 2 k_{i_1} \alpha$.

Hence, $(p_{i_1}-p_{i_0}) + (q_{i_1}-q_{i_0}) \alpha = \frac{\varepsilon 2 \alpha }{y} (k_{i_0}-k_{i_1}) \neq 0$.

So, $\frac{\varepsilon 2 \alpha }{y}$ can be written as $r+s\alpha$ where $r$ and $s$ are rational numbers.

Now, if $n\geq 1$, we still have $(p_{i_n}-p_{i_0}) + (q_{i_n}-q_{i_0}) \alpha = (r+s\alpha) (k_{i_0}-k_{i_n})$.

Because $(1,\alpha)$ is free over $\Q$, we have
\begin{itemize}
\item $(p_{i_n}-p_{i_0}) = r (k_{i_0}-k_{i_n})$
\item $(q_{i_n}-q_{i_0}) = s (k_{i_0}-k_{i_n}) \neq 0$ (remember that $q_n$ is strictly increasing)
\end{itemize}

Thus, by dividing, $$\frac{r}{s} = \frac{p_{i_n}-p_{i_0}}{q_{i_n}-q_{i_0}} = \frac{p_{i_n}}{q_{i_n}} (1-\frac{p_{i_0}}{p_{i_n}}) (\frac{1}{1-\frac{q_{i_0}}{q_{i_n}}}) \cge \alpha \in \R \setminus \Q $$

leading to a contradiction.

For the second case, if $y \in [1,2[$, it is exactly the same
(take the point $A(0,2)$ as the origin and reverse Figure \ref{unfold}).\\

Thus, {\bf the billiard $\P_\alpha$ lacks the finite blocking
property}.

\clearpage

\section{Conclusion}

In \cite{M}, we study Hiemer and Snurnikov's proof: it works
for rational billiards with discrete translation group (such
billiards are called {\em almost integrable}).
Then we generalize the notion of finite blocking property to
translation surfaces (see \cite{MT} for precise definitions). 
With an analogous construction to the one described above, we
obtain the following results:

\begin{thm}
Let $n\geq 3$ be an integer. The following assertions are
equivalent:
\begin{itemize}
\item the regular $n$-gon has the finite blocking property.

\item the right-angled triangle with an angle equal to  $\pi/n$
has the finite blocking property.

\item $n\in\{3,4,6\}$.
\end{itemize}
\end{thm}

\begin{thm}
A translation surface that admits cylinder decomposition of
commensurable moduli in two transversal directions has the finite
blocking property if and only if it is a torus branched covering.
\end{thm}

\begin{cor}
A Veech surface has the finite blocking property if and only if it
is a torus branched covering.
\end{cor}

Note that torus branched coverings are the analogue (in the
vocabulary of translation surfaces) of almost integrable
billiards.

We also provide a local sufficient condition for a translation
surface to fail the finite blocking property: it enables us to
give a complete classification for the L-shaped surfaces and a
density result in the space of translation surfaces in every genus
$g\geq 2$.

\end{document}